\documentclass[12pt]{article}%
\usepackage{graphicx}
\usepackage{amsmath}
\usepackage{amsfonts}
\usepackage{amssymb}%
\newtheorem{theorem}{Theorem}

\newtheorem{lemma}[theorem]{Lemma}

\newenvironment{proof}[1][Proof]{\noindent\textbf{#1.} }{\ \rule{0.5em}{0.5em}}
\begin{document}

\title{\textbf{Self-averaging property of queuing systems.}}
\author{Alexandre Rybko\\Institute for the Information Transmission Problems,\\Russian Academy of Sciences, Moscow, Russia\\rybko@iitp.ru
\and Senya Shlosman\\Centre de Physique Theorique, CNRS,\\Luminy Case 907,\\13288 Marseille, Cedex 9, France\\shlosman@cpt.univ-mrs.fr
\and Alexandre Vladimirov\\Institute for the Information Transmission Problems,\\Russian Academy of Sciences, Moscow, Russia\\vladim@iitp.ru}
\maketitle

\begin{abstract}
We establish the averaging property for a queuing process with one server,
$M(t)/GI/1$. It is a new relation between the output flow rate and the input
flow rate, crucial in the study of the Poisson Hypothesis. Its implications
include the statement that the output flow always possesses more regularity
than the input flow.


\end{abstract}

\section{Introduction}

The Poisson Hypothesis deals with large queuing systems. It is the statement
that for certain large networks the input flow to any given node is
approximately Poissonian with constant rate.

In the paper \cite{RS1} the Poisson Hypothesis is proven for some simple
queuing networks. One of the main technical ingredients of this proof is the
following non-linear averaging relation for $M(t)/GI/1$ queuing process with
one server:
\begin{equation}
b\left(  t\right)  =\left[  \lambda\left(  \cdot\right)  \ast q_{\lambda
,t}\left(  \cdot\right)  \right]  \left(  t\right)  . \label{200}%
\end{equation}
Here $\ast$ stays for convolution: for two functions $a\left(  \cdot\right)
,b\left(  \cdot\right)  $ it is defined as
\[
\left[  a\left(  \cdot\right)  \ast b\left(  \cdot\right)  \right]  \left(
t\right)  =\int a\left(  t-x\right)  b\left(  x\right)  \,dx.
\]
In order to explain the rest of the relation $\left(  \ref{200}\right)  $ we
now introduce notation for our server. Here $\lambda\left(  t\right)  ,$
$-\infty<t<\infty$ is the rate of the Poisson process of moments of arrivals
of customers to our server. If the server is busy, the customer waits in line
for his turn; the service discipline is First-In-First-Out (FIFO). The service
time $\eta_{i}$ for the $i$-th customer is supposed to be random, while the
random variables $\eta_{i}$ are independent identically distributed random
variables, with common distribution $\eta.$ Upon completion of the service the
customer exits the system. The exit flow is of course not Poissonian in
general; yet its rate $b\left(  t\right)  $ is defined. The claim $\left(
\ref{200}\right)  $ is that the functions $\lambda$ and $b$ are related via
the convolution with the kernel $q_{\lambda,t}\left(  x\right)  ,$ which is
supported by positive semiaxis:%
\[
q_{\lambda,t}\left(  x\right)  =0\text{ for }x<0,
\]
and, what is of crucial importance, is stochastic: for every $t$%
\begin{equation}
\int q_{\lambda,t}\left(  x\right)  ~dx=1. \label{201}%
\end{equation}
The function $q_{\lambda,t}\left(  x\right)  $ depends on the function
$\lambda\left(  y\right)  $ only via its restriction to the semiline $\left\{
y\leq t\right\}  .$ (Of course, it depends also on the law of $\eta.$) If the
system is not overloaded, the family $q_{\lambda,t}\left(  \cdot\right)  $ of
probability measures is also compact: for every $\varepsilon>0$ there exists a
threshold $K$ such that $\int_{0}^{K}q_{\lambda,t}\left(  x\right)
~dx>1-\varepsilon,$ uniformly in $t.$ As it is explained in \cite{RS1}, the
averaging relation $\left(  \ref{200}\right)  $ with \textit{stochastic}
kernel $q_{\lambda,t}$ does not hold in general for other disciplines. Some
other examples of self-averaging violation are presented in Sect. \ref{L}. The
importance of $\left(  \ref{200}\right)  $ lies in the fact that it implies
the rate $b$ is in a sense \textquotedblleft smoother\textquotedblright\ than
$\lambda$; in particular, it\ implies that $\sup b\leq\sup\lambda$ and $\inf
b\geq\inf\lambda,$ and it \textquotedblleft almost\textquotedblright\ implies
that the last inequalities are in fact strict.

The proof of $\left(  \ref{200}\right)  -\left(  \ref{201}\right)  ,$ given in
\cite{RS1}, is quite complicated, being based on the validity of a certain
combinatorial identity dealing with rod placements on the real line (see also
\cite{RS2} and \cite{A}). The independence of the service times $\eta_{i}$ is
very important in this proof. The purpose of the present paper is to extend
the relations $\left(  \ref{200}\right)  -\left(  \ref{201}\right)  $ to the
case when the sequence $\eta_{i}$ is not necessarily independent, but has a
weaker property of being a stationary ergodic process. We replace the
combinatorial identity by a stochastic one. In a sense, both in \cite{RS1} and
here we are making use of Fubini theorem, and we use here a different choice
of coordinates, which enables our extension to the dependent case. The
applications of the present result to the Poisson Hypothesis will be presented elsewhere.

It is noteworthy to repeat that the result of the present paper is again based
on an identity, this time a stochastic one. To formulate it, consider the
random variable $V,$ which is a functional of the service process trajectory,
$\omega,$ and which is defined as follows:

\begin{itemize}
\item if the realization $\omega$ is such that the server is idle at the
moment $t=0,$ then $V\left(  \omega\right)  =0;$

\item if the server is occupied at $t=0,$ then let us introduce

$\hat{\eta}\left(  \omega\right)  $ to be the total service time required by
the customer, who is being served at this moment $t=0,$

$\hat{t}\left(  \omega\right)  <0$ to be the moment of the beginning of the
occupation period of the server, which period contains the moment $t=0,$

and finally put
\[
V\left(  \omega\right)  =\frac{1}{\lambda\left(  \hat{t}\left(  \omega\right)
\right)  \hat{\eta}\left(  \omega\right)  }.
\]

\end{itemize}

We claim now that
\begin{equation}
\mathbb{E}\left(  V\left(  \omega\right)  \right)  \equiv1, \label{30}%
\end{equation}
provided only that the server is not overloaded, plus some general technical
conditions. Being very general, the identity $\left(  \ref{30}\right)  $ is
surprisingly non-evident, similarly to the rod placement identity!

\section{The main result}

\subsection{Notation}

Let $\mathcal{P}$ be a Poisson point process on $\mathbb{R}^{1}$ of arrival
moments of the customers. It is a probability measure on the set
$\Omega^{\prime}=\left\{  ...<z_{-1}<z_{0}<z_{1}<...\right\}  $ of
double-infinite sequences $z_{i}\in\mathbb{R}^{1},$ which are locally finite
subsets of $\mathbb{R}^{1}$. Every such Poisson process $\mathcal{P}$ is
defined by the choice of a measure $dm$ on $\mathbb{R}^{1},$ $\mathcal{P}%
=\mathcal{P}_{m},$ and we will suppose that
\[
dm=\lambda\left(  x\right)  dx,
\]
where $\lambda>0$ is the rate of the process $\mathcal{P}_{m}=\mathcal{P}%
_{\lambda}.$

{\small Strictly speaking, the Poisson process is a measure on locally finite
subsets }$\phi\subset\mathbb{R}^{1}${\small . We consider it as a measure on
sequences by defining the point }$z_{0}$ {\small to be the smallest positive
point in }$\phi.$

We further assume that the process $\mathcal{T}$ of \textit{positive }service
times $\left\{  \ldots,\eta_{-1},\eta_{0},\eta_{1},\ldots\right\}  ,$
independent of $\mathcal{P}_{\lambda},$ is given. We will assume that
$\mathcal{T}$ is stationary and ergodic.

The total process we thus are interested in, is the process $\mathcal{S=P}%
_{\lambda}\times\mathcal{T},$ which is a probability measure on the set
$\Omega=\left\{  \omega=\ldots,\left(  z_{-1},\eta_{-1}\right)  ,\left(
z_{0},\eta_{0}\right)  ,\left(  z_{1},\eta_{1}\right)  ,\ldots\right\}  .$

\subsection{Nonlinear shift}

In the special case when the rate $\lambda$ equals to the constant $\ell,$ the
process $\mathcal{S}$ is ergodic with respect to the shift transformation
$T_{t}$ on $\Omega:$%
\begin{align*}
&  T_{t}\left(  \ldots,\left(  z_{-1},\eta_{-1}\right)  ,\left(  z_{0}%
,\eta_{0}\right)  ,\left(  z_{1},\eta_{1}\right)  ,\ldots\right) \\
&  =\ldots,\left(  z_{-1}+t,\eta_{-1}\right)  ,\left(  z_{0}+t,\eta
_{0}\right)  ,\left(  z_{1}+t,\eta_{1}\right)  ,\ldots.
\end{align*}
In the case of $\lambda$ arbitrary the same is true once $T_{t}$ is replaced
by the non-linear shift $\theta_{t}$. It is defined as follows: for every
$x,t\in\mathbb{R}^{1}$ we define $\theta_{t}\left(  x\right)  \in
\mathbb{R}^{1}$ as the only value for which%
\[
\int_{x}^{\theta_{t}\left(  x\right)  }\lambda\left(  x\right)  dx=t.
\]
Clearly,%
\[
\frac{d}{dt}\theta_{t}\left(  x\right)  \Bigm|_{t=0}=\frac{1}{\lambda\left(
x\right)  },
\]
so the shift $\theta_{t}\left(  x\right)  $ is the same as traveling along the
vector field $\frac{dx}{\lambda\left(  x\right)  }$ from the location $x$ for
the time duration $t.$ We suppose that the following integrals diverge:
\begin{equation}
\int_{-\infty}^{0}\lambda\left(  x\right)  dx=\infty,\ \ \int_{0}^{\infty
}\lambda\left(  x\right)  dx=\infty, \label{34}%
\end{equation}
then the measure $\lambda\left(  x\right)  dx$ is invariant under every
transformation $\theta_{t}.$ The claims that the process $\mathcal{S}$ is
invariant and ergodic under the transformation%
\begin{align*}
&  \theta_{t}\left(  \ldots,\left(  z_{-1},\eta_{-1}\right)  ,\left(
z_{0},\eta_{0}\right)  ,\left(  z_{1},\eta_{1}\right)  ,\ldots\right) \\
&  =\ldots,\left(  \theta_{t}\left(  z_{-1}\right)  ,\eta_{-1}\right)
,\left(  \theta_{t}\left(  z_{0}\right)  ,\eta_{0}\right)  ,\left(  \theta
_{t}\left(  z_{1}\right)  ,\eta_{1}\right)  ,\ldots.
\end{align*}
are immediate.

\subsection{The exit flow}

Let $\omega=\left(  \ldots,\left(  z_{-1},\eta_{-1}\right)  ,\left(
z_{0},\eta_{0}\right)  ,\left(  z_{1},\eta_{1}\right)  ,\ldots\right)
\in\Omega.$ In case that for all $i$ we have \textit{no conflicts:}%
\begin{equation}
z_{i}\geq z_{i-k}+\eta_{i-k}\text{ for all }k>0, \label{11}%
\end{equation}
we define the exit moments $E\left(  \omega\right)  =\left\{  \ldots
<y_{-1}<y_{0}<y_{1}<\ldots\right\}  $ in the evident way by%
\begin{equation}
y_{i}=z_{i}+\eta_{i}. \label{12}%
\end{equation}
Otherwise we need to \textit{resolve the conflicts}. To do so we first
introduce the set $I\left(  \omega\right)  \subset\mathbb{Z}^{1}$ of all
indices $i,$ for which the relation $\left(  \ref{11}\right)  $ holds. Assume
$\omega$ is such that the set $I\left(  \omega\right)  $ is double-infinite
sequence $\left\{  \ldots<i_{-1}<i_{0}<i_{1}<\ldots\right\}  .$ We define the
sequence $R\omega=\left(  \ldots,\left(  Rz_{-1},\eta_{-1}\right)  ,\left(
Rz_{0},\eta_{0}\right)  ,\left(  Rz_{1},\eta_{1}\right)  ,\ldots\right)  $ in
the following way: if $j\in I\left(  \omega\right)  ,$ then $Rz_{j}=z_{j}.$
For other $j$-s we have $i_{k}<j<i_{k+1}$ for some $k,$ and we put%
\[
Rz_{j}=z_{i_{k}}+\eta_{i_{k}}+\ldots+\eta_{j-1}%
\]
(Lindley equation). In case $R\omega$ has no conflicts, we can again use
$\left(  \ref{12}\right)  .$ Otherwise, if the set $I\left(  R\omega\right)  $
is again double-infinite, we can define the configuration $R^{2}\omega,$ and
so on. Assume that the configuration $\omega$ is such that

\begin{enumerate}
\item the configurations $R^{n}\omega$ are defined for all $n\geq1,$

\item for every $j$ the sequence $R^{n}z_{j}$ stabilizes at some finite
$n=n\left(  \omega,j\right)  .$ We denote by $\bar{R}$ the limiting
transformation, which will be called the \textit{conflict resolution
operator.}
\end{enumerate}

\noindent Then we define the exits $E\left(  \omega\right)  $ by%
\[
y_{i}=\bar{R}z_{i}+\eta_{i}.
\]
Note that if the configuration $\omega\in\Omega$ violates only finitely many
of the following sequence of conditions:%
\begin{align}
\eta_{-1}+z_{-1}  &  <0,\ \ \eta_{-1}+\eta_{-2}+z_{-2}<0,...,\nonumber\\
\sum_{i=-k}^{-1}\eta_{i}+z_{-k}  &  <0,..., \label{36}%
\end{align}

\noindent then the sequence $R^{n}\left(  \omega\right)  $ stabilizes
pointwise, at every location, so the operator $\bar{R}\left(  \omega\right)  $
is defined. We denote by $\tilde{\Omega}\subset\Omega$ the subset of
configurations thus defined, and we put $\bar{\Omega}=\cap_{-\infty<t<\infty
}\theta_{t}\tilde{\Omega}.$

Our main assumption on the process $\mathcal{S}$ is the condition%
\begin{equation}
\mathcal{S}\left(  \bar{\Omega}\right)  =\mathcal{S}\left(  \cap
_{-\infty<t<\infty}\theta_{t}\tilde{\Omega}\right)  =1. \label{141}%
\end{equation}
It can be viewed as a natural generalization of the usual condition the
process not to be overloaded. In what follows we will call the condition
$\left(  \ref{141}\right)  $ the \textit{no-overload }condition. Note that the
no-overload condition implies that the probability of the server to be idle at
any given time moment is positive.

The main example of the no-overloaded process $\mathcal{S}$ is obtained by
imposing the following restriction on the rate of the Poisson process
$\lambda:$
\begin{equation}
\ell\equiv\limsup_{T\rightarrow\infty}\frac{1}{T}\int_{-T}^{0}\lambda\left(
x\right)  \ dx<1, \label{13}%
\end{equation}
while taking $\mathcal{T}$ to be stationary and shift-ergodic, with
\begin{equation}
\mathbb{E}\left(  \eta_{i}\right)  \equiv1. \label{140}%
\end{equation}
The proof of that statement is the content of the Lemma \ref{idle} below. In a
way, it shows that the no-overload condition is dynamically insensitive, in
the sense of \cite{BHPS}.

In what follows, every point $z_{i}\in\omega$, such that $R^{n}z_{i}=z_{i}$
for all $n,$ will be called the \textit{point of the beginning of the cluster}
of $\omega$, or the \textit{head of the cluster}. We will call the set
\[
\mathcal{C}\left(  z_{i},\omega\right)  =\left\{  z_{i}+\eta_{i},z_{i}%
+\eta_{i}+\eta_{i+1},\ldots,z_{i}+\eta_{i}+\eta_{i+1}\ldots+\eta
_{j-1}\right\}
\]
the cluster of $z_{i},$ where $j>i$ is the index of the next after $z_{i}$
point $z_{j}$ of the beginning of a cluster. (In case $R\omega=\omega,$ all
the points $z_{i}$ are heads of clusters, while each cluster contains just one
point, $z_{i}+\eta_{i}.$) The segment $\left[  z_{i},z_{i}+\eta_{i}+\eta
_{i+1}\ldots+\eta_{j-1}\right]  $ will be called the \textit{support }of the
cluster $\mathcal{C}\left(  z_{i},\omega\right)  ,$ while the segment $\left[
z_{i}+\eta_{i}+\eta_{i+1}\ldots+\eta_{j-1},z_{j}\right]  $ -- the gap between
the consecutive clusters.

\subsection{The Averaging Theorem}

Let $b\left(  t\right)  $ be the rate of the process $E\left(  \omega\right)
:$%
\[
b\left(  t\right)  =\lim_{\Delta t\rightarrow0}\frac{1}{\Delta t}\mathbf{\Pr
}\left\{  E\left(  \omega\right)  \cap\left[  t,t+\Delta t\right]
\neq\varnothing\right\}  .
\]
The kernels $q_{\lambda,t}\left(  x\right)  $ are the same kernels which were
used in \cite{RS1}. They are defined as follows. Let $e\left(  u\right)  $ be
the probability that our server is idle at the time $u.$ (Note that the
dependence of $e\left(  u\right)  $ on $\lambda$ is only via $\left\{
\lambda\left(  x\right)  ,x\leq u\right\}  .$) Now define the function
$c\left(  u,x\right)  $ as follows. Let us condition on the event that the
server is idle just before time $u,$ while at $u$ the customer arrives. Under
this condition define
\begin{equation}
c\left(  u,x\right)  =\lim_{h\searrow0}\frac{1}{h}\mathbf{\Pr}\left\{
\begin{array}
[c]{c}%
\text{the server is never idle during }\left[  u,u+x\right]  ;\text{ }\\
\text{during }\left[  u+x,u+x+h\right]  \text{ the server gets }\\
\text{through with some customer}%
\end{array}
\right\}  . \label{008}%
\end{equation}
Then
\begin{equation}
q_{\lambda,t}\left(  x\right)  =e\left(  t-x\right)  c\left(  t-x,x\right)  .
\label{005}%
\end{equation}

\begin{theorem}
\label{1} Let the service time process $\mathcal{T}=\left\{  \ldots,\eta
_{-1},\eta_{0},\eta_{1},\ldots\right\}  $ be stationary and ergodic, while the
Poisson process $\mathcal{P}_{\lambda}$ is defined by the \textbf{continuous
positive} rate function $\lambda,$ such that the no-overload assumption
$\left(  \ref{141}\right)  $ and condition $\left(  \ref{34}\right)  $ hold.
Then for the kernels $q_{\lambda,t}$ we have%
\begin{equation}
b\left(  t\right)  =\left[  \lambda\left(  \cdot\right)  \ast q_{\lambda
,t}\left(  \cdot\right)  \right]  \left(  t\right)  , \label{18}%
\end{equation}
$-\infty<t<\infty.$ The kernels $q_{\lambda,t}$ depend on the function
$\lambda$ only via restrictions $\lambda\Bigm|_{(-\infty,t]}.$ Moreover, they
are stochastic: for every $t$%
\begin{equation}
\int_{0}^{\infty}q_{\lambda,t}\left(  x\right)  ~dx=1, \label{19}%
\end{equation}
while $q_{\lambda,t}\left(  x\right)  =0$ for $x<0.$
\end{theorem}

For the case of the process $\mathcal{T}$ to be that of independent
identically distributed random variables, this theorem was proven in sections
5 and 6 of \cite{RS1}. As in that paper, the main problem is to show the
relation $\left(  \ref{19}\right)  .$ The combinatorial counting, applied
there, is valid only in the independent case.

\section{Proof of the Theorem}

\subsection{The representation for the kernels $q_{\lambda,t}$}

In this section we will obtain another representation of the kernels
$q_{\ast,\ast},$ which will elucidate more clearly its stochastic nature.

Let us denote by $I_{x}$ the indicator of the event $H\left(  x,\Delta
x\right)  \subset\Omega,$ that the intersection $E\left(  \omega\right)
\cap\left[  x,x+\Delta x\right]  \neq\varnothing.$ Then
\[
b\left(  x\right)  =\lim_{\Delta x\rightarrow0}\frac{1}{\Delta x}\int_{\Omega
}I_{x}\left(  \omega\right)  d\mathcal{S}\left(  \omega\right)  .
\]
By shift-invariance of $\mathcal{S}$ we also can write
\begin{equation}
b\left(  x\right)  =\lim_{\Delta x\rightarrow0}\frac{1}{\Delta x}\int_{\Omega
}\left(  \int_{0}^{1}I_{x}\left(  \theta_{t}\omega\right)  dt\right)
d\mathcal{S}\left(  \omega\right)  .\label{16}%
\end{equation}

Let us fix any $\omega$ and consider all the moments $t\in\left[  0,1\right]
,$ for which $E\left(  \theta_{t}\omega\right)  \cap\left[  x,x+\Delta
x\right]  \neq\varnothing.$ We will call them the \textit{hitting }moments.
The set of all hits will be denoted by $\tau\left(  \omega\right)
\subset\left[  0,1\right]  .$ Without loss of generality we can assume that
for every $t\in\tau\left(  \omega\right)  $ the intersection $E\left(
\theta_{t}\omega\right)  \cap\left[  x,x+\Delta x\right]  $ consists of just
one point. (This is certainly the case, once $\Delta x$ is chosen to be small enough.)

Evidently, the set $\tau\left(  \omega\right)  \subset\left[  0,1\right]  $ of
hitting moments is a union of disjoint segments, $\tau\left(  \omega\right)
=\cup_{r=1}^{s\left(  \omega\right)  }D_{r}$. Clearly,%
\[
\int_{0}^{1}I_{x}\left(  \theta_{t}\omega\right)  dt=\sum_{r=1}^{s\left(
\omega\right)  }l\left(  D_{r}\right)  ,
\]
where $l$ stays for the length of the segment. For every $t\in\tau\left(
\omega\right)  $ we define now the index $i\left(  t\right)  \in\mathbb{Z}%
^{1},$ to be the one satisfying the relation%
\[
\mathcal{C}\left(  \theta_{t}z_{i\left(  t\right)  },\theta_{t}\omega\right)
\cap\left[  x,x+\Delta x\right]  \neq\varnothing,
\]
while the index $j\left(  t\right)  $ is the one for which
\[
y_{j\left(  t\right)  }\left(  \theta_{t}\omega\right)  \in\left[  x,x+\Delta
x\right]  .
\]
These indices are well-defined. In words, the point $z_{i\left(  t\right)
}\in\omega$ is the one which becomes the head of the hitting cluster after the
shift $\theta_{t}$ is applied, while the point $y_{j\left(  t\right)  }$ is
just the intersection $\mathcal{C}\left(  \theta_{t}z_{i\left(  t\right)
},\theta_{t}\omega\right)  \cap\left[  x,x+\Delta x\right]  $. Let us further
partition the set $\tau\left(  \omega\right)  \subset\left[  0,1\right]  $
into maximal segments of constancy of the function $i\left(  t\right)  .$ Let
us denote this partition by $\tau\left(  \omega\right)  =\cup_{k=1}^{u\left(
\omega\right)  }C_{k},$ while $i\left(  k\right)  $ will denote the (constant)
value of the function $i\left(  t\right)  $ when $t\in C_{k}.$ Evidently, we
have%
\[
\int_{0}^{1}I_{x}\left(  \theta_{t}\omega\right)  dt=\sum_{k=1}^{u\left(
\omega\right)  }l\left(  C_{k}\right)  .
\]
In general, the partition $\cup_{k=1}^{u\left(  \omega\right)  }C_{k}$ is
finer than the partition $\cup_{r=1}^{s\left(  \omega\right)  }D_{r}.$
However, once $\Delta x$ is small enough, they are the same, provided that we
know a priori that the point $y_{j\left(  t\right)  }\left(  \theta_{t}%
\omega\right)  $ moves continuously with time. In what follows we are assuming
this continuity, and we postpone the proof of it till the end of the present subsection.

Let us compute the length $l\left(  C_{k}\right)  $ of the segment
$C_{k}=\left[  c_{k},e_{k}\right]  \subset\left[  0,1\right]  $. At the moment
$t=c_{k}$ the cluster $\mathcal{C}\left(  \theta_{t}z_{i\left(  k\right)
},\theta_{t}\omega\right)  $ starts to hit the segment $\left[  x,x+\Delta
x\right]  ,$ which means in particular that
\begin{equation}
\theta_{c_{k}}z_{i\left(  k\right)  }+\eta_{i\left(  k\right)  }%
+\eta_{i\left(  k\right)  +1}\ldots+\eta_{j_{k}}=x \label{20}%
\end{equation}
for some appropriate $j_{k}\geq i\left(  k\right)  ,$ $j_{k}=j_{k}\left(
\omega\right)  .$ As the time $t$ increases from $c_{k}$ to $e_{k},$ the point
$\theta_{t}z_{i\left(  k\right)  }+\eta_{i\left(  k\right)  }+\eta_{i\left(
k\right)  +1}\ldots+\eta_{j_{k}}$ moves from the initial value $x$ up to the
terminal value $x+\Delta x.$ Let us compute the time $l\left(  C_{k}\right)  $
it takes. By definition, the point $\theta_{t}z_{i\left(  k\right)  }%
+\eta_{i\left(  k\right)  }+\eta_{i\left(  k\right)  +1}\ldots+\eta_{j_{k}}$
moves with the velocity
\[
\frac{1}{\lambda\left(  \theta_{t}z_{i\left(  k\right)  }\right)  }.
\]
Therefore, we have%
\[
\int_{c_{k}}^{e_{k}}\frac{dt}{\lambda\left(  \theta_{t}z_{i\left(  k\right)
}\right)  }=\Delta x,
\]
so%
\[
l\left(  C_{k}\right)  =e_{k}-c_{k}=\lambda\left(  \theta_{\tilde{t}_{k}%
}z_{i\left(  k\right)  }\right)  \Delta x
\]
for some $\tilde{t}_{k}\in$ $C_{k},$ due to the Mean Value Theorem. Hence,
from $\left(  \ref{16}\right)  $ we have
\begin{equation}
b\left(  x\right)  =\int_{\Omega}\left(  \sum_{k=1}^{\mathfrak{S}\left(
1,\omega,x\right)  }\lambda\left(  \theta_{t_{k}\left(  x\right)  }z_{i\left(
k\right)  }\right)  \right)  d\mathcal{S}\left(  \omega\right)  ., \label{17}%
\end{equation}
where the times $t_{k}\left(  x\right)  ,$ $1\leq k\leq\mathfrak{S}\left(
1,\omega,x\right)  $ are all the moments in$\ \left[  0,1\right]  ,$ at which
some cluster, $\mathcal{C}\left(  \theta_{t_{k}\left(  x\right)  }z_{i\left(
k\right)  },\theta_{t_{k}\left(  x\right)  }\omega\right)  ,$ contains the
point $x$. In case $\mathfrak{S}\left(  1,\omega,x\right)  $ vanishes, we
define the sum $\sum_{k=1}^{0}$ to be zero. Now we see that the claim $\left(
\ref{18}\right)  $ of our Theorem holds, with the kernel
\begin{equation}
q_{\lambda,x}=\int_{\Omega}\left(  \sum_{k=1}^{\mathfrak{S}\left(
1,\omega,x\right)  }\delta_{\eta_{i\left(  k\right)  }+\eta_{i\left(
k\right)  +1}\ldots+\eta_{j_{k}}}\right)  d\mathcal{S}\left(  \omega\right)  ,
\label{21}%
\end{equation}
see the relation $\left(  \ref{20}\right)  .$ (Compare also with the similar
relations (34) and (40) from \cite{RS1}.)

{\small A little thought shows that indeed the r.h.s. of }$\left(
{\small \ref{21}}\right)  $ {\small depends strongly on }$\lambda
\Bigm|_{(-\infty,x]}.$ {\small For example, if the rate }$\lambda$ {\small is
small in a segment }$\left[  x_{0},x\right]  ,${\small suitably long,}
{\small then the sum of }$\eta${\small -s in the subscript of the
delta-function will typically have only one summand.}

We conclude this subsection by proving the continuity statement used above.

\begin{lemma}
\label{cont} The exit moments $y_{i}\left(  \theta_{t}\omega\right)  $ of the
configuration $\theta_{t}\omega,$ $i=0,\pm1,...$ -- are continuous functions
of $t,$ once $\omega$ is taken from $\bar{\Omega}.$
\end{lemma}

\begin{proof}
Let us consider only the case $i=0,$ and suppose that $y_{0}\left(
\theta_{t=0}\omega\right)  =1,$ say.

Clearly, if we would have impose the restriction that $\lambda\geq c>0,$ then
our claim is immediate, since every point would move with a speed not
exceeding $c^{-1}.$ However we know only that $\lambda>0,$ so our points can
have arbitrarily high speeds, and in principle it is feasible that clusters
successively accelerate each other and produce an infinite speed somewhere. As
we will show, that does not happen once $\omega\in\bar{\Omega}.$

Let $\mathcal{C}\left(  z_{i\left(  t\right)  },\theta_{t}\omega\right)  $ be
the cluster containing the exit $y_{0}\left(  \theta_{t}\omega\right)  .$
Evidently we will be done once we know that for any $T$ there exists a
constant $K\left(  T,\omega\right)  ,$ such $i\left(  t\right)  \geq K\left(
T,\omega\right)  $ for all $t\in\left[  0,T\right]  .$ Indeed, that means that
the movement of the point $y_{0}\left(  \theta_{t}\omega\right)  $ is
determined only by finitely many other points, while all of them have finite
speeds. To see the existence of $K\left(  T,\omega\right)  $ we note that for
any $\omega\in\bar{\Omega},$ we can write that%
\[
y_{0}\left(  \omega\right)  =\sup_{-\infty<k\leq0}\left(  z_{k}+\eta
_{k}+...+\eta_{0}\right)  .
\]
Moreover, we know that
\begin{equation}
\limsup_{k\rightarrow-\infty}\left(  z_{k}+\eta_{k}+...+\eta_{0}\right)
=-\infty. \label{40}%
\end{equation}
In particular, for the location $y_{0}\left(  \theta_{T}\omega\right)
>y_{0}\left(  \theta_{0}\omega\right)  =1$ we have for some finite $k\left(
T\right)  $ that
\[
\theta_{T}z_{k\left(  T\right)  }+\eta_{k\left(  T\right)  }+...+\eta
_{0}=y_{0}\left(  \theta_{T}\omega\right)  .
\]
Because of $\left(  \ref{40}\right)  $ we know that for some other (negative)
value $K\left(  T\right)  <k\left(  T\right)  $ and for all $k\leq K\left(
T\right)  $
\[
\theta_{T}z_{k}+\eta_{k}+...+\eta_{0}<1.
\]
But then we have that for all $t\leq T$%
\[
\theta_{t}z_{k}+\eta_{k}+...+\eta_{0}<1
\]
since $\theta_{t}z_{k}$ is increasing in $t.$ That estimate means that no
point $\theta_{t}z_{k}$ with $k\leq K\left(  T\right)  $ can be in the same
cluster with the point $y_{0}\left(  \theta_{t}\omega\right)  \left(
\geq1\right)  $ for all $t\in\left[  0,T\right]  .$
\end{proof}

\subsection{Counting of exit moments}

To prove $\left(  \ref{19}\right)  ,$ it remains to establish the following

\begin{theorem}
Suppose the process $\mathcal{S}$ satisfies the no-overload property $\left(
\ref{141}\right)  .$ Then for every $x$%
\begin{equation}
\mathbb{E}\left(  \mathfrak{S}\left(  1,\omega,x\right)  \right)  =1.
\label{33}%
\end{equation}

\end{theorem}

\begin{proof}
For $T\geq0$ let us introduce the random variables $\mathfrak{S}%
(T,\omega,x)\ $as the number of indices $i$ such that $y_{i}(\omega)<x$ and
$y_{i}(\theta_{T}\omega)>x.$ For $T\leq0$ define similarly the random
variables $\mathfrak{S}^{\prime}(T,\omega,x)$ to be the number of indices $i$
such that $y_{i}(\omega)>x$ and $y_{i}(\theta_{T}\omega)<x.$ Evidently,
$\mathfrak{S}(T,\omega,x)=\mathfrak{S}^{\prime}(-T,\theta_{T}\omega,x).$ By
shift-invariance of $\mathcal{S},$ we have $\mathbb{E}\left(  \mathfrak{S}%
^{\prime}(-T,\theta_{T}\omega,x)\right)  =\mathbb{E}\left(  \mathfrak{S}%
^{\prime}(-T,\omega,x)\right)  .$ So for our purposes it is enough to show
that
\begin{equation}
\mathbb{E}\left(  \mathfrak{S}\left(  1,\omega,x\right)  \right)  \geq1
\label{101}%
\end{equation}
and
\begin{equation}
\mathbb{E}\left(  \mathfrak{S}^{\prime}(-1,\omega,x)\right)  \leq1.
\label{102}%
\end{equation}

If instead of termination points $y$-s, crossing the location $x$ under a
time-shift, we will count the number of $z$-points, crossing it, we obtain
similarly the random variables $\mathfrak{R}(T,\omega,x),$ defined for
$T\geq0,$ and $\mathfrak{R}^{\prime}(T,\omega,x),$ defined for $T\leq0.$ It is
immediate from our definition of the non-linear shift, that
\begin{equation}
\mathbb{E}\left(  \mathfrak{R}\left(  1,\omega,x\right)  \right)
=\mathbb{E}\left(  \mathfrak{R}^{\prime}(-1,\omega,x)\right)  =1. \label{105}%
\end{equation}

By ergodic theorem (see e.g. \cite{CFS}), for $\mathcal{S}$-a.e. $\omega$%
\begin{equation}
\mathbb{E}\left(  \mathfrak{S}\left(  1,\omega,x\right)  \right)
=\lim_{T\rightarrow\infty}\frac{1}{T}\mathfrak{S}(T,\omega,x), \label{103}%
\end{equation}%
\begin{equation}
\mathbb{E}\left(  \mathfrak{R}\left(  1,\omega,x\right)  \right)
=\lim_{T\rightarrow\infty}\frac{1}{T}\mathfrak{R}(T,\omega,x). \label{104}%
\end{equation}
Let us define the queue length $\mathfrak{Q}(\omega,x)$ as the number of
indices $i$ such that $z_{i}(\omega)\leq x$ and $y_{i}(\omega)>x$. It follows
from the above definitions that for any $T$%
\[
\mathfrak{S}(T,\omega,x)\geq\mathfrak{R}(T,\omega,x)-\mathfrak{Q}(\omega,x).
\]
Therefore the relations $\left(  \ref{105}\right)  -\left(  \ref{104}\right)
$ imply $\left(  \ref{101}\right)  .$

In the same way we have, that%
\[
\mathbb{E}\left(  \mathfrak{S}^{\prime}\left(  -1,\omega,x\right)  \right)
=\lim_{T\rightarrow\infty}\frac{1}{T}\mathfrak{S}^{\prime}(-T,\omega,x),
\]%
\[
\mathbb{E}\left(  \mathfrak{R}^{\prime}\left(  -1,\omega,x\right)  \right)
=\lim_{T\rightarrow\infty}\frac{1}{T}\mathfrak{R}^{\prime}(-T,\omega,x),
\]
and%
\[
\mathfrak{R}^{\prime}(-T,\omega,x)\geq\mathfrak{S}^{\prime}(-T,\omega
,x)-\mathfrak{Q}(\omega,x),
\]
which together imply $\left(  \ref{102}\right)  .$
\end{proof}

\subsection{Time integral approximation for $\mathfrak{S}$}

In this section we obtain under the conditions of the Theorem \ref{1} the
infinitesimal version of the identity $\left(  \ref{33}\right)  $ -- the
identity $\left(  \ref{30}\right)  .$ To do it, we \textit{approximate} the
function $\mathfrak{S}\left(  T,\omega,x\right)  $ by a time integral
average,
\begin{equation}
\frac{1}{T}\int_{0}^{T}V_{x}\left(  \theta_{t}\left(  \omega\right)  \right)
~dt, \label{24}%
\end{equation}
and we claim that the following function $V_{x}\left(  \omega\right)  $ is suitable:

$i)$ if no cluster of $\omega$ has the point $x$ inside its support, then
$V_{x}\left(  \omega\right)  =0,$

$ii)$ in the opposite case we have
\[
z\left(  \omega,x\right)  +\eta_{i\left(  \omega,x\right)  }\ldots+\eta
_{j-1}<x<z\left(  \omega,x\right)  +\eta_{i\left(  \omega,x\right)  }%
\ldots+\eta_{j}%
\]
for some cluster $\mathcal{C}\left(  z\left(  \omega,x\right)  ,\omega\right)
$ of $\omega$ and some $j=j\left(  \omega,x\right)  \geq i\left(
\omega,x\right)  ;$ we take
\[
V_{x}\left(  \omega\right)  =\frac{1}{\lambda\left(  z\left(  \omega,x\right)
\right)  \eta_{j\left(  \omega,x\right)  }}.
\]
To explain the relation between the integral $\left(  \ref{24}\right)  $ and
the number of summands in $\left(  \ref{21}\right)  ,$ let $t^{\prime}$
$\ $(resp., $t^{\prime\prime}$) be the moment when the above rod $\eta_{j}$
starts (resp., ends) to cover the point $x,$ i.e.%
\[
z\left(  \theta_{t^{\prime}}\left(  \omega\right)  ,x\right)  +\eta_{i\left(
\omega,x\right)  }\ldots+\eta_{j}=x,\text{ resp. }z\left(  \theta
_{t^{\prime\prime}}\left(  \omega\right)  ,x\right)  +\eta_{i\left(
\omega,x\right)  }\ldots+\eta_{j-1}=x.
\]
At the moment $t\in\left(  t^{\prime},t^{\prime\prime}\right)  $ the point $x$
moves relative to the rod $\eta_{j}$ with velocity $\left(  \lambda\left(
z\left(  \theta_{t}\left(  \omega\right)  ,x\right)  \right)  \right)  ^{-1},$
hence%
\[
\int_{t^{\prime}}^{t^{\prime\prime}}\frac{dt}{\lambda\left(  z\left(
\theta_{t}\left(  \omega\right)  ,x\right)  \right)  }=\eta_{j\left(
\omega,x\right)  }.
\]
Therefore%
\[
\int_{t^{\prime}}^{t^{\prime\prime}}V_{x}\left(  \theta_{t}\left(
\omega\right)  \right)  ~dt=1,
\]
and so%
\begin{equation}
\left\vert \int_{0}^{T}V_{x}\left(  \theta_{t}\left(  \omega\right)  \right)
~dt-\mathfrak{S}\left(  T,\omega,x\right)  \right\vert \leq2, \label{32}%
\end{equation}
where the difference is due to the influence of the rods at the ends of the
interval $\left[  0,T\right]  ,$ one per end.

Therefore the expectation%
\begin{equation}
\mathbb{E}\left(  \frac{1}{T}\int_{0}^{T}V_{x}\left(  \theta_{t}\left(
\omega\right)  \right)  ~dt\right)  \rightarrow1\text{ as }T\rightarrow\infty,
\label{31}%
\end{equation}
because of $\left(  \ref{33}\right)  .$ (As we will see soon, the expectation
$\mathbb{E}\left(  \frac{1}{T}\int_{0}^{T}V_{x}\left(  \theta_{t}\left(
\omega\right)  \right)  ~dt\right)  $ in fact equals to $1$ for every $T>0.$)
On the other hand, due to the ergodic theorem,%
\[
\mathbb{E}\left(  \frac{1}{T}\int_{0}^{T}V_{x}\left(  \theta_{t}\left(
\omega\right)  \right)  ~dt\right)  =\lim_{\Upsilon\rightarrow\infty}\frac
{1}{\Upsilon}\int_{0}^{\Upsilon}\left(  \frac{1}{T}\int_{0}^{T}V_{x}\left(
\theta_{t+\tau}\left(  \omega\right)  \right)  ~dt\right)  ~d\tau,
\]
for $\mathcal{S}$-almost every $\omega.$ But, evidently, the r.h.s. limit does
not depend on $T,$ and moreover
\[
\lim_{\Upsilon\rightarrow\infty}\frac{1}{\Upsilon}\int_{0}^{\Upsilon}\left(
\frac{1}{T}\int_{0}^{T}V_{x}\left(  \theta_{t+\tau}\left(  \omega\right)
\right)  ~dt\right)  ~d\tau=\lim_{\Upsilon\rightarrow\infty}\frac{1}{\Upsilon
}\int_{0}^{\Upsilon}V_{x}\left(  \theta_{\tau}\left(  \omega\right)  \right)
~d\tau.
\]
Therefore for every $T>0$%
\[
\mathbb{E}\left(  \frac{1}{T}\int_{0}^{T}V_{x}\left(  \theta_{t}\left(
\omega\right)  \right)  ~dt\right)  =\lim_{\Upsilon\rightarrow\infty}\frac
{1}{\Upsilon}\int_{0}^{\Upsilon}V_{x}\left(  \theta_{\tau}\left(
\omega\right)  \right)  ~d\tau=1,
\]
due to $\left(  \ref{31}\right)  .$ That proves the identity $\left(
\ref{30}\right)  .\medskip$

We conclude this section by presenting the proof of the

\begin{lemma}
\label{idle} Suppose the rate $\lambda$ of the Poisson process $\mathcal{P}%
_{\lambda}$ satisfies
\[
\ell\equiv\limsup_{T\rightarrow\infty}\frac{1}{T}\int_{-T}^{0}\lambda\left(
x\right)  \ dx<1,
\]
while the process $\mathcal{T}$ is stationary and shift-ergodic, with
\[
\mathbb{E}\left(  \eta_{i}\right)  \equiv1.
\]
Then the set $\bar{\Omega}\subset\Omega$ of trajectories $\omega,$ for which
the conflict resolution operator $\bar{R}\left(  \theta_{t}\omega\right)
\ $is defined for all $t,$ has full measure:%
\[
\mathcal{S}\left(  \bar{\Omega}\right)  =1.
\]

\end{lemma}

This result should be compared with a classical result of Loynes, \cite{L},
see also the books \cite{B} and \cite{SD}.

\begin{proof}
It is easy to see that the domain $\tilde{\Omega}\subset\Omega$ (see relation
$\left(  \ref{141}\right)  $) has full measure. Indeed, the validity of all
except finitely many of the relations $\left(  \ref{36}\right)  $ holds
$\mathcal{S}$-almost surely, due to the Strong Law of Large Numbers. This law
holds here since $\mathbb{E}\left(  \eta\right)  =1,$ while $\ell<1.$
Moreover, the subset $\hat{\Omega}\subset\tilde{\Omega}$ of configurations
$\omega,$ satisfying the relation
\begin{equation}
\limsup_{k\rightarrow\infty}\left(  \sum_{i=-k}^{-1}\eta_{i}+z_{-k}\right)
=-\infty, \label{37}%
\end{equation}
also has full measure, for the same reason.

Next, let us show that the intersection $\bar{\Omega}=\cap_{-\infty<t<\infty
}\theta_{t}\tilde{\Omega}$ also has full measure. Clearly, the countable
intersection $\cap_{t\in\mathbb{Z}^{1}}\theta_{t}\hat{\Omega}$ has measure
one. So we will be done once we show that for every $T$%
\[
\cap_{t\in\mathbb{Z}^{1}}\theta_{t}\hat{\Omega}\subset\theta_{T}\tilde{\Omega
}.
\]
This is the same as to claim that for any $\omega\in\cap_{t\in\mathbb{Z}^{1}%
}\theta_{t}\hat{\Omega}$ and any $T$ we have $\theta_{-T}\omega\in
\tilde{\Omega}.$ This will be established once we show a stronger statement,
that for any $t>0,$ any $\omega=\left\{  \left(  z_{i},\eta_{i}\right)
\right\}  \in\hat{\Omega}$ we have $\theta_{-t}\omega\in\hat{\Omega}.$ To
check this inclusion we have to consider the $k\rightarrow\infty$ asymptotics
of the sums%
\[
\sum_{i=-k}^{n\left(  \omega,t\right)  }\eta_{i}+\theta_{-t}z_{-k},
\]
where $n\left(  \omega,t\right)  $ is the largest index $i,$ satisfying the
relation $\theta_{-t}z_{i}<0.$ But since $\theta_{-t}z_{-k}\leq z_{-k},$ we
have evidently that%
\[
\limsup_{k\rightarrow\infty}\left(  \sum_{i=-k}^{n\left(  \omega,t\right)
}\eta_{i}+\theta_{-t}z_{-k}\right)  =-\infty
\]
as well.
\end{proof}

\section{\label{L} Self-averaging not always holds}

The example of the service discipline without self-averaging property, given
in \cite{RS1}, is very simple. All the customers arriving between $n$
\textit{A.M. }and $n+1$ \textit{A.M. }are leaving the server at $n+1$
\textit{A.M. }sharp. Self-averaging violation is evident in this case.

However, the above discipline may not look very natural. So below we present
an example which look nicer. It has the property that whenever the queue is
non-empty the server is busy, and that, once served, the customer leaves the server.

The example is the following. Suppose the server gets two kinds of clients:
\textquotedblleft slow\textquotedblright\ and \textquotedblleft
fast\textquotedblright. The slow one needs time $L$ for its service, while the
fast one needs time $l\ll L,$ both times are non-random. The probability that
a given client turns out to be slow is $\frac{1}{2},$ say, and the sequence of
service times is iid (with two values). The server has a preference to gather
fast clients into big groups and to work on the group as a whole. It is doing
so by having two modes -- slow and fast. This means the following: in the slow
mode, whenever the queue contains both the fast and the slow clients, the
server picks a slow one, provided the number of the fast ones does not exceed
a certain value, $F.$ At the moment when the number of waiting fast clients
becomes $F,$ the server switches into fast mode. It starts serving the fast
clients and works on them till none is left in the queue, and then switches
back to slow mode and resumes its service of the slow customer.

To see that the self-averaging is violated, let us take the rate function
$\lambda\left(  t\right)  $ of the input Poisson flow to be zero for $t<0$ and
$\Lambda\gg1$ for $t\geq0.$ We are going to explain that if $L,F$ and
$\Lambda$ are large enough, while $l$ is small enough, then there exists a
time moment $T>0,$ at which the exit flow rate $b\left(  T\right)  $ exceeds
the value $\Lambda,$ thus violating any possibility of self-averaging.

To see it, let us consider the random moment $\tau_{1}$ of the beginning of
the service of the first slow customer. Its distribution is well localized
around $t=0;$ in fact, it tends to $\delta_{0}$ as $\Lambda\rightarrow\infty.$
The time $\tau_{2}$ it takes to wait till $F$ fast customers arrive after the
moment $\tau_{1},$ is such that $\mathbb{E}\left(  \tau_{2}\right)  =\frac
{2F}{\Lambda}.$ If $L$ is big enough -- namely, if $L\geq L_{0}\left(
F\right)  $ for some $L_{0}$ -- then it is very likely that during the time
interval $\left[  \tau_{1},\tau_{1}+\tau_{2}\right]  $ the server is never
idle and is serving only slow clients. That means that at the moment $\tau
_{1}+\tau_{2}$ the service of $F$ waiting fast clients will start. Then during
the time interval $\left[  \tau_{1}+\tau_{2},\tau_{1}+\tau_{2}+Fl\right]  $ it
will happen $F$ exits (of fast clients). Due to the Law of Large Numbers, this
(random) time interval belongs to the interval $\left[  \frac{2F}{\Lambda
}\left(  1-\varepsilon\right)  ,\frac{2F}{\Lambda}\left(  1+\varepsilon
\right)  +Fl\right]  $ with very high probability, provided $F$ is large.
Therefore the exit rate has to exceed the value
\[
\frac{F}{\frac{4F}{\Lambda}\varepsilon+Fl}=\frac{\Lambda}{4\varepsilon+\Lambda
l}%
\]
somewhere inside this interval, which in turn exceeds $\Lambda$ for $l$ and
$\varepsilon$ small enough.

The above server has some kind of a memory. We conclude this paper by
formulating a conjecture about what we call memoryless disciplines. The server
will be called memoryless, if its strategy at any given moment can depend only on

\begin{itemize}
\item the number of the clients waiting in the queue at this moment,

\item the times the clients need,

\item the order at which the clients were arriving to the server.
\end{itemize}

\noindent In particular, we exclude the server to have its own clock. In
addition, the server can not be idle if there are clients waiting in the
queue. Finally, if some client has its service started, then it goes on
without interruptions and delays to its end.

We believe that if the service discipline is memoryless, then the
self-averaging does hold for it, and it can be proven by the extension of the
methods of the present paper.

\textbf{Acknowledgment. }\textit{We would like to thank our colleagues -- in
particular, F. Kelly, K. Khanin, Yu. Peres, S. Pirogov, O. Zeitouni -- for
valuable discussions and remarks, concerning this paper. A. Rybko would like
to acknowledge the generous financial support of the Leverhulme Trust. He is
also grateful to the Centre de Physique Theorique in Luminy, for the financial
support during his visit in the Fall of 2004, when this work was initiated. S.
Shlosman would like to acknowledge the generous financial support of the
Department of Mathematics of UC Berkeley and the uplifting atmosphere of the
Mathematical Sciences Research Institute during the Program "Probability,
Algorithms \& Statistical Physics", in the Spring, 2005, when part of that
work was done.}

\end{document}